\documentclass[12pt,twoside]{article}
\usepackage{amsxtra,amssymb,amsthm,amsmath,latexsym}

\theoremstyle{plain}

\def\oH{\buildrel\circ\over H}
\def\oH1{\buildrel\circ\over H\kern-.02in{}^1}
\def\oH1{\buildrel\circ\over H\kern-.02in{}^1}

\def\d{\delta}
\def\ep{\epsilon}

\begin{document}
\title{A new discrepancy principle
\thanks{Math subject classification: 47H15, 45G10,35B25}
\thanks{key words:  ill-posed problems, discrepancy principle }                  
}

\author{                          
A.G. Ramm\\       
 Mathematics Department, Kansas State University, \\
 Manhattan, KS 66506-2602, USA\\
ramm@math.ksu.edu\\}

\date{}

\maketitle\thispagestyle{empty}

1. {\bf Inroduction}

 The aim of this note is to  prove a {\it new discrepancy
principle}. The advantage of the new discrepancy principle 
compared with the known one consists of solving a minimization
problem approximately, rather than exactly, and in the proof of 
a stability result.
To explain this in 
more detail, let us  recall the usual  discrepancy principle, 
which can be stated as follows.
Consider an operator eqution 
$$Au=f,
 \eqno{(1)}
$$ 
where $A:H\to H$
 is a bounded linear operator on a Hilbert space $H$, and assume 
that the 
range $R(A)$ is not closed, so that problem (1) is 
ill-posed.
Assume that $f=Ay$ where $y$ is the minimal-norm solution to (1),
and that  noisy data $f_\d$ are given, such that $||f_\d -f||\leq 
\d$. One wants to construct a stable approximation to $y$, given 
 $f_\d$.
The variational regularization method for solving this problem 
consists of solving the minimization problem
$$
F(u):=||Au-f_\d||^2 +\ep||u||^2=\min. 
 \eqno{(2)}
$$
It is well known that problem (2) has a solution and this solution is 
unique (see e.g. [1]). Let $u_{\d, \ep}$ solve (2). Consider the  equation
for finding $\ep=\ep(\d)$:
$$ ||Au_{\d, \ep}-f_\d||=C\d,
 \eqno{(3)}
$$
where $C=const>1$. Equation (3) is the usual discrepancy 
principle. One can prove that equation (3) determines 
 $\ep=\ep(\d)$  uniquely, $\ep(\d)\to 0$ as $\d\to 0$,
and $u_\d:=u_{\d, \ep(\d)}\to y$ as $\d\to 0$. This justifies 
the usual discrepancy principle for choosing the regularization 
parameter (see [1] and [2] for various  justifications
of this principle, [3] for the dynamical systems method
for stable solution of equation (1), and [4] for a method
of solving nonlinear ill-posed problems).

The drawback of this principle consists of the necessity to
solve problem (2) exactly. The other 
drawback is the lack of information concerning stability of the 
solution to (3): if one solves (2) approximately in some sense, 
will
the element $u_{\d, \ep(\d)}$ ( with $\ep(\d)$ being an 
approximate solution to (3)) converge to $y$?

Our aim is to formulate and justify a new  discrepancy principle
which deals with the both issues mentioned above.

 Our basic result is:

{\bf Theorem 1.} {\it Assume: 

i) $A$ is a bounded linear 
operator in a Hilbert space $H$, 

ii) equation $Au=f$ is 
solvable, and
$y$ is its minimal-norm solution, 

iii) $||f_\d -f||\leq \d$,
$||f_\d||>C\d$, where $C>1$ is a constant.

Then: 

j) equation (3)
is solvable for $\ep$ for any fixed $\d>0$, where
 $u_{\d,\ep}$ is any element satisfying inequality 
$F(u_{\d,\ep})\leq m+(C^2-1-b)\d^2$,   
$F(u):=||A(u)-f_\d||^2+\ep 
||u||^2,$ $m=m(\d,\ep):=inf_{u} F(u)$,  $b=const>0$, and $C^2>1+b,$

and

jj) if $\ep=\ep(\d)$ solves (3), and 
$u_\d:=u_{\d,\ep(\d)}$, then $\lim_{\d \to 0}||u_\d-y||=0$.}

In Section 2 proof of Theorem 1 is given.

2. {\bf Proof}

{\bf Proof of Theorem 1. } 
Let us first prove the existence of a solution to (3). 
We {\it claim} that the function $h(\d, \ep):=||Au_{\d,\ep}-f_\d||$ is
greater 
than $C\d$ 
for sufficiently large $\ep$, and smaller than $C\d$ for sufficiently 
small $\ep$.
If this is proved, then the 
continuity of $h(\d, \ep)$ with respect to $\ep$ on $(0,\infty)$
implies that the equation 
$h(\d, \ep)=C\d$ has a solution.

Let us prove the {\it claim}. As $\ep\to \infty$, we use the inequality:
$$\ep ||u_{\d,\ep}||^2\leq F(u_{\d,\ep})\leq m+(C^2-1-b)\d^2 \leq 
F(0)+(C^2-1-b)\d^2,$$
and, as $\ep\to 0$, we use another inequality: 
$$||Au_{\d,\ep}-f_\d||^2<F(u_{\d,\ep})\leq m+(C^2-1-b)\d^2 \leq 
F(y)+(C^2-1-b)\d^2= 
\ep ||y||^2+(C^2-b)\d^2.$$ 
This inequality implies
$$h^2(\d,\ep)< \ep ||y||^2+(C^2-b)\d^2.$$
 As $\ep \to \infty$, one gets
$||u_{\d,\ep}||\leq \frac c {\sqrt {\ep}}\to 0,$ where $c>0$ is a constant
depending on $\d$. Thus, by the continuity of $A$, one obtains
$$\lim_{\ep\to \infty}h(\d, \ep)= ||A(0)-f_\d||=||f_\d||>C\d.$$ 

As $\ep \to 0$, one gets
$$\liminf_{\ep \to 0}h(\d,\ep)=\liminf_{\ep \to 0}( \ep 
||y||^2+(C^2-b)\d^2)^{1/2}<C\d.$$ 
Therefore equation $h(\d,\ep)=C\d$ has a solution $\ep=\ep(\d)>0.$

Let us now prove that if $u_\d:=u_{\d,\ep(\d)}$, then 
$$\lim_{\d\to 0}||u_\d-y||=0. 
 \eqno{(4)}
$$
From the estimate 
$$F(u_\d=||Au_\d-f_\d||^2+\ep||u_\d||^2\leq C^2\d^2+\ep||y||^2,
$$ 
and from  (3), it follows that 
$$||u_\d||\leq ||w||,
 \eqno{(5)}
$$ 
where $w$ is {\it any} solution to (1). We will use (5) with $w=y$ and
$w=U$,
where $U$ is a solution to (1) constructed below, and $y$ is a minimal-norm
solution to (1).

Thus, one may assume that $u_\d\rightharpoonup U$,
and from (3) it follows that 
$Au_\d\to f$ as $\d\to 0$.
 
This implies, as we prove below (see (8)), that
$$AU=f. 
 \eqno{(6)}
$$ 
We also prove below that from (5)
it follows that 
$$\lim_{\d\to 0}||u_\d-U||=0, \quad ||U||\leq ||y||.
 \eqno{(7)}
$$ 
The minimal 
norm solution to  equation (1) is unique. Consequently, (7) implies 
$U=y$. Thus, (4) holds.

Let us now prove (6). We have
$$(f,v)=\lim_{\d\to 0}(Au_\d,v)=\lim_{\d\to 0}(u_\d, A^*v)=(U,
A^*v)=(AU,v)\quad
\forall v\in H.
 \eqno{(8)}
$$
Since $v$ is arbitrary, (8) implies (6).

Finally, we prove (7). We have $u_\d\rightharpoonup U$. Thus,
$||U||\leq \lim\inf_{\d\to
0}||u_\d||$. Inequality (5) implies $ \lim\sup_{\d\to
0}||u_\d||\leq ||U||$.
Consequently, $||u_\d||\to  ||U||$. It is well known that the weak
convergence
together with convergence of the norms imply in a Hilbert space strong
convergence.
Therefore $\lim_{\d\to 0}||u_\d-U||=0$. Taking $w=y$ in (5), and then
passing to the limit $\d\to 0$ in (5),  yields
inequality (7). Thus, both parts of (7) are established.
 
Since $U$ solves equation (1), and 
$||U||\leq ||y||$, it follows that $U=y$, and  (4) holds.

Theorem 1 is proved $\Box$.

\end{document}